\documentclass[12pt] {article}
\usepackage{amsmath,amsthm}
\usepackage{amssymb,latexsym}
\usepackage{amscd}
\usepackage{url}
\title{On modules over valuations.}
\date{}
\author{Semyon Alesker \footnote{Partially supported by ISF grant 701/08.}
\\  { \normalsize Department of Mathematics, Tel Aviv University, Ramat Aviv}
\\  { \normalsize 69978 Tel Aviv, Israel }
\\ {\normalsize e-mail: semyon@post.tau.ac.il}}

\def\alp{\alpha}

\def\to{\rightarrow}
\def\qed { Q.E.D. }

\def\CC{\mathbb{C}}

\swapnumbers
\newtheorem{theorem}{Theorem}[section]
\newtheorem{corollary}[theorem]{Corollary}
\newtheorem{lemma}[theorem]{Lemma}
\newtheorem{proposition}[theorem]{Proposition}
\newtheorem{claim}[theorem]{Claim}
\theoremstyle{definition}

\newtheorem{remark}[theorem]{Remark}

\theoremstyle{proposition-definition}
\newtheorem{proposition-definition}[theorem]{Proposition-Definition}

\def\cf{{\cal F}}

\def\ca{{\cal A}} \def\cb{{\cal B}} 
 \def\ce{{\cal E}} \def\cf{{\cal F}}
 \def\ch{{\cal H}} 
 \def\ck{{\cal K}} \def\cl{{\cal L}}
\def\cm{{\cal M}}  \def\co{{\cal O}}

\def\cv{{\cal V}} \def\cw{{\cal W}} 
 
\def\inj{\hookrightarrow }
\def\surj{\twoheadrightarrow}
\def\vi{V^\infty}

\def\shv{\cv_X^\infty}
\begin{document}
\maketitle

\def\vxo{\cv_{x_0}}

\begin{abstract}
To any smooth manifold $X$ an algebra of smooth valuations
$V^\infty(X)$ was associated in
\cite{part1}-\cite{part4},\cite{part3}. In this note we initiate a
study of $V^\infty(X)$-modules. More specifically we study finitely
generated projective modules in analogy to the study of vector
bundles on a manifold. In particular it is shown that for a compact
manifold $X$ there exists a canonical isomorphism between the
$K$-ring constructed out of finitely generated projective
$V^\infty(X)$-modules and the classical topological $K^0$-ring
constructed out of vector bundles over $X$.
\end{abstract}

\section{Introduction.}\label{S:introduction}
Let $X$ be a smooth manifold of dimension $n$.\footnote{All
manifolds are assumed to be countable at infinity, i.e. presentable
as a union of countably many compact subsets. In particular they are
paracompact.} In \cite{part1}-\cite{part4},\cite{part3} the notion
of a smooth valuation on $X$ was introduced. Roughly put, a smooth
valuation is a $\CC$-valued finitely additive measure on compact
submanifolds of $X$ with corners, which satisfies in addition some
extra conditions. We omit here the precise description of the
conditions due to their technical nature. Let us notice that basic
examples of smooth valuations include any smooth measure on $X$ and
the Euler characteristic. There are many other natural examples of
valuations coming from convexity, integral, and differential
geometry. We refer to recent lecture notes \cite{alesker-barcelona},
\cite{fu-barcelona}, \cite{bernig-AIG} for an overview of the
subject, examples, and applications.

The space $V^\infty(X)$ of all smooth valuations is a Fr\'echet
space. It has a canonical product making $V^\infty(X)$ a commutative
associative algebra over $\CC$ with a unit element (which is the
Euler characteristic).

\hfill

In this note we initiate a study of modules over $V^\infty(X)$. Our
starting point is the analogy to the following well known fact due
to Serre and Swan \cite{serre}, \cite{swan}: if $X$ is compact, then
the category of smooth vector bundles of finite rank over $X$ is
equivalent to the category of finitely generated projective modules
over the algebra $C^\infty(X)$ of smooth functions (the functor in
one direction is given by taking global smooth sections of a vector
bundle).

In order to state our main results we need to remind a few general
facts about valuations on manifolds. We have a canonical
homomorphism of algebras
\begin{eqnarray}\label{E:hom-alg}
V^\infty(X)\to C^\infty(X)
\end{eqnarray}
given by the evaluation on points, i.e. $\phi\mapsto [x\mapsto
\phi(\{x\})]$. This is an epimorphism. The kernel, denoted by $W_1$,
is a nilpotent ideal of $V^\infty(X)$:
$$(W_1)^{n+1}=0.$$

Next, smooth valuations form a sheaf of algebras which is denoted by
$\cv^\infty_X$: for an open subset $U\subset X$,
$$\cv^\infty_X(U)=V^\infty(U),$$
where the restriction maps are obvious. We denote by $\co_X$ the
sheaf of $C^\infty$-smooth functions on $X$. Then the map
(\ref{E:hom-alg}) gives rise to the epimorphism of sheaves
\begin{eqnarray}\label{E:hom-sh}
\cv^\infty_X\surj \co_X.
\end{eqnarray}

Recall now the notion of a projective module. Let $A$ be a
commutative associative algebra with a unit. An $A$-module $M$ is
called {\itshape projective} if $M$ is a direct summand of a free
$A$-module, i.e. there exists an $A$-module $N$ such that $M\oplus
N$ is a free $A$-module (not necessarily of finite rank). It is easy
to see that if $M$ is in addition finitely generated then $M$ is a
direct summand of a free $A$-module of finite rank.

Let $\ca$ be a sheaf of algebras on a topological space $X$. A sheaf
$\cm$ of $\ca$-modules is called a {\itshape locally projective}
$\ca$-module if any point $x\in X$ has an open neighborhood $U$ such
that $\cm(U)$ is a projective $\ca(U)$-module.

Let us denote by $Proj_f\cv_X^\infty-mod$ the full subcategory of
$\cv^\infty_X$-modules  consisting of locally projective
$\cv^\infty_X$-modules of finite rank. Let us denote by
$Proj_f\vi(X)-mod$ the full subcategory of the category of
$V^\infty(X)$-modules consisting of projective $\vi(X)$-modules of
finite rank. In Section \ref{S:Loc_Free} we prove the following
result.
\begin{theorem}\label{T:1}
Let $X$ be a smooth manifold.

(1) Any locally projective $\cv^\infty_X$-module of finite rank is
locally free.

(2) Assume in addition that $X$ is compact. Let $\ce$ be a locally
free $\cv^\infty_X$-module of finite rank. Then there exists another
locally free $\cv^\infty_X$-module $\ch$ of finite rank such that
$\ce\oplus \ch$ is isomorphic to $(\cv_X^\infty)^N$ for some natural
number $N$.

(3) Assume again that $X$ is compact. Then the functor of global
sections
$$\Gamma\colon Proj_f\cv_X^\infty-mod\to Proj_f V^\infty(X)-mod$$
is an equivalence of categories.
\end{theorem}
Notice that all the statements of the theorem are completely
analogous to the classical situation of vector bundles (whose spaces
of sections are projective finitely generated
$C^\infty(X)$-modules). For example a version of (2) for vector
bundles says that any vector bundle is a direct summand of a free
bundle. A classical version of (3) is the above mentioned theorem of
Serre-Swan. The method of proof of Theorem \ref{T:1} is a minor
modification of the proof for the analogous statement for vector
bundles.

\hfill

To formulate our next main result observe that to any
$\cv^\infty_X$-module we can associate an $\co_X$-module via
\begin{eqnarray}\label{E:val-bund}
\cm\mapsto \cm\otimes_{\cv^\infty_X}\co_X,
\end{eqnarray}
where $\co_X$ is considered as $\cv^\infty_X$-module via the
epimorphism (\ref{E:hom-sh}). Clearly under this correspondence
locally free $\cv^\infty_X$-modules of finite rank are mapped to
locally free $\co_X$-modules of equal rank, i.e. to vector bundles.
\begin{theorem}\label{T:2}
Assume that $X$ is a compact manifold. Let $N$ be a natural number.
The map (\ref{E:val-bund}) induces a bijection between the
isomorphism classes of locally free $\cv^\infty_X$-modules of rank
$N$ and isomorphism classes of vector bundles of rank $N$.
\end{theorem}
Theorem \ref{T:2} is proved in Section \ref{S:Vect_Bun}. The proof
is an application of general results of Grothendieck \cite{groth} on
non-abelian cohomology of topological spaces and the existence of a
finite decreasing filtration on $\cv^\infty_X$ such that the
associated graded sheaf is a sheaf of $\co_X$-modules.

\hfill

{\bf Acknowledgement.} I thank M. Borovoi for useful discussions on
non-abelian cohomology, and F. Schuster for numerous remarks on the
first version of the paper.

\section{Locally free sheaves over valuations.}\label{S:Loc_Free}

A sheaf $\ce$ of $\cv^\infty_X$-modules is called locally projective
of finite rank if every point $x\in X$ has a neighborhood $U$ such
that there exists a sheaf of $\cv^\infty_U$-modules $\cf$ with the
property that $\ce|_U\oplus \cf$ is isomorphic to $(\cv^\infty_X)^N$
for some natural number $N$.

The technique used in the proofs of most of the results of this
section is rather standard and is a simple modification of that from
\cite{atiyah-lectures}.

\begin{proposition}\label{P:bun0}
On a manifold $X$ any locally projective $\shv$-module of finite
rank is locally free.
\end{proposition}
{\bf Proof.} Fix a point $x_0\in X$. Let us denote for brevity
$\vxo:=\cv_{X,x_0}^\infty$ (resp. $\co_{X,x_0}$) the stalk at $x_0$
of the sheaf $\shv$ (resp. $\co_X$). Let $\ce$ be a locally
projective $\shv$-module of finite rank. Consider its stalk
$\ce_{x_0}$ as $\vxo$-module. Then there exists a $\vxo$-module
$\cf$ such that
$$\ce_{x_0}\oplus \cf\simeq \vxo^N$$
for some natural number $N$. Consider the idempotent endomorphism of
the $\vxo$-module
$$e\colon \vxo^N\to \vxo^N$$
given by the projection onto $\ce_{x_0}$. Thus $e^2=e$. Notice that
$\vxo$ is a local ring with the maximal ideal
$$m:=\{\phi\in\vxo\, |\, \phi(\{x_0\})=0\}.$$
Clearly $\vxo/m=\CC$.

We have
\begin{eqnarray}\label{bun1}
\CC^N=\vxo^N\otimes_{\vxo}(\vxo/m)=(\ce_{x_0}\otimes_{\vxo}\vxo/m)\oplus(\cf\otimes_{\vxo}
\vxo/m).
\end{eqnarray}

Let us choose a basis
$$\xi_1',\dots,\xi_k',f_1',\dots,f_{N-k}'$$
of $\CC^N$ such that the $\xi_i'$'s form a basis of the first
summand in the right hand side of (\ref{bun1}), and the $f_j'$'s
form a basis of the second summand. Let $\tilde \xi_i\in \vxo^N,\,
\tilde f_j\in \vxo^N$ be their lifts. Define finally
\begin{eqnarray*}
\xi_i:=e(\tilde\xi_i)\in \ce_{x_0},\\
f_j:=(1-e)(\tilde f_j)\in \cf.
\end{eqnarray*}

It is clear that
\begin{eqnarray*}
\xi_i\equiv \tilde\xi_i \,\, mod(m),\\
f_j\equiv \tilde f_j \, \, mod(m).
\end{eqnarray*}

Consider the morphism of $\cv_{x_0}$-modules $\theta\colon \vxo^k\to
\ce_{x_0}$ given by
$$\theta(\phi_1,\dots,\phi_k):=\sum_{i=1}^k\phi_i\xi_i.$$
Consider also another morphism $\tau\colon \vxo^{N-K}\to \cf$ given
by
$$\tau(\psi_1,\dots,\psi_{N-k})=\sum_{j=1}^{N-k}\psi_jf_j.$$
Now define a morphism of $\vxo$-modules by
$$\sigma:=\theta\oplus \tau\colon \vxo^N=\vxo^k\oplus \vxo^{N-k}\to
\ce_{x_0}\oplus \cf\simeq \vxo^N.$$ We claim that $\sigma$ is an
isomorphism. It is equivalent to the property that $\det(\sigma)\in
\vxo$ is invertible. In order to see this it suffices to show that
$(\det\sigma)(\{x_0\})\ne 0$. But the last condition is satisfied
since
$$\sigma\otimes Id_{\vxo/m}\colon \vxo^N\otimes_{\vxo}\vxo/m\to
\vxo^N\otimes_{\vxo}\vxo/m$$ is an isomorphism $\CC^N\to \CC^N$
since by construction
$$\xi_1'=\xi_1(\{x_0\}),\dots,\xi_k'=\xi_k(\{x_0\}),f_1'=f_1(\{x_0\}),\dots,f_{N-k}'=f_{N-K}(\{x_0\})$$
form a basis of $\CC^N$.

Since $\sigma\colon \cv_{x_0}^N\to\cv_{x_0}^N$ is an isomorphism, it
follows that there exists an open neighborhood $U$ of $x_0$ and an
isomorphism of $\cv_U^\infty$-modules
$$\tilde\sigma\colon (\cv_U^\infty)^N\to (\cv_U^\infty)^N$$
which extends $\sigma$, i.e. $\sigma$ is the stalk of $\tilde\sigma$
at $x_0$. It follows that $\theta$ extends to an isomorphism
$$\tilde\theta\colon (\cv_U^\infty)^k\tilde\to \ce|_U$$
of $\cv_U^\infty$-modules (and similarly for $\tau$). \qed

\hfill

\begin{lemma}\label{bun2}
Let $\ce$ be a locally free $\cv^\infty_X$-module of finite rank
$N$. Let $\xi_1,\dots,\xi_k\in H^0(X,\ce)$ be chosen such that for
every point $x_0\in X$ their images $\bar\xi_1,\dots,\bar\xi_k$ in
$\ce\otimes_{\cv^\infty_X}\cv^\infty_X/m_{x_0}\simeq \CC^N$ form a
linearly independent sequence. Consider the morphism of
$\cv^\infty_X$-modules
$$f\colon (\cv^\infty_X)^k\to \ce,$$
given by $f(\phi_1,\dots,\phi_k)=\sum_{i=1}^k\phi_i\xi_i$.

Then $f\colon (\cv^\infty_X)^k\to Im(f)$ is an isomorphism, and
$\ce/Im(f)$ is a locally free $\cv^\infty_X$-module of rank $N-k$.
\end{lemma}
{\bf Proof.} Let us denote for brevity $\cv:=\cv^\infty_X$. The
statement is local on $X$. Fix $x_0\in X$. We can choose
$\eta_1,\dots,\eta_{N-k}\in H^0(\ce)$ such that their images
$$\bar\xi_1,\dots,\bar\xi_k,\bar\eta_1,\dots,\bar\eta_{N-k}\in
\ce\otimes_{\cv}\cv/m_{x_0}$$ form a basis. Consider a morphism of
$\cv$-modules $g\colon \cv^k\oplus \cv^{N-k}\to \ce$ given by
$$g(\phi_1,\dots,\phi_k;\psi_1,\dots,\psi_{N-k})=\sum_{i=1}^k\phi_i\xi_i+\sum_{j=1}^{N-k}\psi_j\eta_j.$$
Clearly $g|_{\cv^k}\equiv f$. In a neighborhood of $x_0$ we may and
will identify $\ce\simeq \cv^N$. Then
$$g\colon \cv^k\oplus \cv^{N-K}=\cv^N\to \cv^N.$$ It is easy to see
that the map
$$g\otimes Id_{\cv/m_{x_0}}\colon \CC^N\simeq (\cv/m_{x_0})^N\to
\CC^N\simeq (\cv/m_{x_0})^N$$ is an isomorphism. Hence $(\det
g)(\{x_0\})\ne 0$. It follows that $\det g\in \cv$ is invertible in
a neighborhood of $x_0$. Hence $g$ is an isomorphism in a
neighborhood of $x_0$. This implies the lemma immediately. \qed

\hfill

\begin{lemma}\label{bun3}
Let $P$ be a locally free $\cv^\infty_X$-module of finite rank. Then
for any $\cv^\infty_X$-module $A$,
$$Ext^i_{\cv^\infty_X-mod}(P,A)=0 \mbox{ for } i>0.$$
\end{lemma}
{\bf Proof.} We abbreviate again $\cv:=\cv^\infty_X$. First notice
that in the category of $\cv$-modules the following two functors
$$F,G\colon \cv-mod\to Vect$$
are naturally isomorphic:
\begin{eqnarray*}
F(A)=Hom_{\cv-mod}(P,A),\\
G(A)=H^0(X,P^*\otimes_\cv A),
\end{eqnarray*}
where $P^*:=\underline{Hom}_{\cv-mod}(P,\cv)$ is the inner Hom as
usual. Indeed the natural morphism
\begin{eqnarray*}
P^*\otimes_\cv A=\underline{Hom}_{\cv-mod}(P,\cv)\otimes_\cv A\to
\underline{Hom}_{\cv-mod}(P,A)
\end{eqnarray*}
is an isomorphism of sheaves. Taking global sections, we get an
isomorphism
$$H^0(X,P^*\otimes_\cv A)\tilde\to
H^0(X,\underline{Hom}_{\cv-mod}(P,A)).$$ But the last space is equal
to $Hom_{\cv-mod}(P,A)$ (see \cite{hartshorne-book}, Ch. II, \S1,
Exc. 1.15).

Consequently $F$ and $G$ have isomorphic derived functors. Hence
$$Ext^i_{\cv-mod}(P,A)\simeq H^i(X, P^*\otimes_\cv A).$$ But the last
group vanishes for $i>0$ by \cite{part4}, Lemma 5.1.2. \qed

\hfill

\begin{corollary}\label{bun4}
Let $$0\to A\to B\to C\to 0$$ be a short exact sequence of
$\cv_X^\infty$-modules. If $C$ is locally free of finite rank, then
this exact sequence splits.
\end{corollary}
{\bf Proof.} Indeed $Ext^1_{\cv-mod}(C,A)=0$ by Lemma \ref{bun3}.
\qed

\hfill

\begin{proposition}\label{bun5}
Let $X$ be a compact manifold. Let $\ce$ be a locally free
$\cv^\infty_X$-module of finite rank. Then there exists another
locally free $\cv^\infty_X$-module $\ch$ of finite rank such that
$$\ce\oplus \ch\simeq (\cv^\infty_X)^N.$$
\end{proposition}
{\bf Proof.} Let us choose a finite open covering $X=\cup_\alp
U_\alp$ such that the sheaf $\ce|_{U_\alp}$ is free for each $\alp$.
Let $\{\phi_\alp\}$ be a partition of unity in the algebra of
valuations subordinate to this covering (it exist by \cite{part4},
Proposition 6.2.1). We can find a finite dimensional subspace
$L_\alp\subset H^0(U_\alp,\ce)$ which generates $\ce|_{U_\alp}$ as
$\cv^\infty_{U_\alp}$-module. Consider $\phi_\alp\cdot L_\alp\subset
H^0(X,\ce)$ (where all sections are extended by zero outside of
$U_\alp$). Then the finite dimensional subspace
$$L:=\sum_\alp\phi_\alp\cdot L_\alp\subset H^0(X,\ce)$$
generates $\ce$ as $\cv$-module (indeed at every $x\in X$ there
exists an $\alp$ such that $\phi_\alp$ is invertible in a
neighborhood of $x$).

Let us choose a basis $\xi_1,\dots,\xi_s$ of $L$. Consider the
morphism of $\cv$-modules $F\colon \cv^s\to \ce$ given by
$$F(\phi_1,\dots,\phi_s)=\sum_{i=1}^s\phi_i\xi_i.$$
Clearly $F$ is an epimorphism of $\cv$-modules. Let $\ca:=Ker(F)$.
By Corollary \ref{bun4} the short exact sequence
$$0\to \ca\to \cv^s\to \ce\to 0$$ splits. Thus $\ce\oplus \ca\simeq
\cv^s$. Hence $\ca$ is locally projective. Hence $\ca$ is locally
free by Proposition \ref{P:bun0}. \qed

\hfill

Let us denote by $Proj_f\cv_X^\infty-mod$ (or just $Proj_f \cv-mod$)
the full subcategory of $\cv-mod$ consisting of locally free
$\cv$-modules of finite rank. Let us denote by $Proj_f\vi(X)-mod$
the category of projective $\vi(X)$-modules of finite rank.
\begin{theorem}\label{P:bun6}
Let $X$ be a compact manifold. Then the functor of global sections
$$\Gamma\colon Proj_f\cv_X^\infty-mod\to Proj_f\vi(X)-mod$$
is an equivalence of categories.
\end{theorem}
{\bf Proof.} We denote again by $\cv:=\cv_X^\infty$. Let $\ca,\cb\in
Proj_f\cv-mod$. First let us show that
$$Hom_{\cv-mod}(\ca,\cb)=Hom_{\vi(X)-mod}(\Gamma(\ca),
\Gamma(\cb)).$$ Both $Hom$ functors respect finite direct sums with
respect to both arguments. Since $\ca,\cb$ are direct summands of
free $\cv$-modules by Proposition \ref{bun5} we may assume that
$\ca=\cb=\cv$. But clearly
\begin{eqnarray*}
Hom_{\cv-mod}(\cv,\cv)=V^\infty(X),\\
Hom_{V^\infty(X)}(\vi(X),\vi(X))=\vi(X).
\end{eqnarray*}
Thus $\Gamma$ is fully faithful.

Let us define a functor in the opposite direction (the localization
functor),
$$G\colon Proj_f\vi(X)-mod\to Proj_f\cv-mod,$$ by
$G(A):=A\otimes_{\vi(X)}\cv$. $G$ is also fully faithful: it
commutes with direct sums, and for trivial $\vi(X)$-modules the
statement is obvious.

The functors $F\circ G$ and $G\circ F$ are naturally isomorphic to
the identity functors. \qed

\def\glv{\underline{GL_N(\cv)}}
\def\glo{\underline{GL_N(\co_X)}}
\section{Isomorphism classes of bundles over
valuations.}\label{S:Vect_Bun} Recall that the sheaf of smooth
valuations $\cv^\infty_X$, which we will denote for brevity by
$\cv$, has a canonical filtration by subsheaves,
$$\cv=\cw_0\supset \cw_1\supset \dots\supset \cw_n.$$
This filtration is compatible with the product, and $\cv/\cw_1\simeq
\co_X$ canonically \cite{part4}. Let us fix a natural number $N$.
Let us denote by $\glv$ (resp. $\glo$) the sheaf on $X$ of
invertible $N\times N$ matrices with entries in $\cv$ (resp.
$\co_X$). We have a natural homomorphism of sheaves of groups
\begin{eqnarray}\label{E:homo-sg-gr}
\glv\to\glo.
\end{eqnarray}

It is well known that isomorphism classes of usual vector bundles
are in bijective correspondence with the (Cech) cohomology set
$H^1(X,\glo)$. Similarly it is clear that locally free $\cv$-modules
of rank $N$ are in bijective correspondence with the set
$H^1(X,\glv)$. The main result of this section is
\begin{theorem}\label{T:biject}
Let $X$ be a compact manifold. The natural map
$$H^1(X,\glv)\to H^1(X,\glo)$$
induced by (\ref{E:homo-sg-gr}) is a bijection. Thus, if $X$ is
compact, the isomorphism classes of rank $N$ locally free
$\cv$-modules are in natural bijective correspondence with
isomorphism classes of rank $N$ vector bundles.
\end{theorem}

\hfill

We will need some preparations before the proof of the theorem.
First we observe that $\glv$ has a natural filtration by subsheaves
of normal subgroups
$$\glv=:\ck_0\supset \ck_1\supset \dots \supset \ck_n,$$
where for any $i>0$, $\ck_i(U):=\{\xi\in \glv(U)|\,\, \xi\equiv I
\mbox{ mod } \cw_i(U)\}$ for any open subset $U\subset X$. We have
the canonical isomorphisms of sheaves of groups:
\begin{eqnarray}
\ck_0/\ck_1\simeq \underline{GL_N(\co_X)},\\\label{IsomorSh-1}
\ck_i/\ck_{i+1}\simeq (\cw_i/ \cw_{i+1})^N \mbox{ for } i>0.
\end{eqnarray}

\def\ung{\underline{G}}
\def\unf{\underline{F}}
\def\unh{\underline{H}}
We have to remind some general results due to Grothendieck
\cite{groth}. Let $\ung$ be a sheaf of groups (not necessarily
abelian) on a topological space $X$. Let $\unf\subset \ung$ be a
subsheaf of {\itshape normal} subgroups. Let $\unh:=\ung/\unf$ be
the quotient sheaf, which is also a sheaf of groups. Notice that the
sheaf of groups $\ung$ acts on $\unf$ by conjugations.

Let $E'$ be a $\ung$-torsor. Let $c':=[E']\in H^1(X,\ung)$ be its
class. Define a new sheaf $\unf(E')$ to be the sheaf associated to
the presheaf
$$U\mapsto \unf(U)\times_{\ung(U)}E'(U).$$
Since $\ung$ acts on $\unf$ by automorphisms, it follows that
$\unf(E')$ is a sheaf of groups.

Grothendieck \cite{groth} has constructed a map
$$i_1\colon H^1(X,\unf(E'))\to H^1(X,\ung),$$
and he has shown (see Corollary after Proposition 5.6.2 in
\cite{groth}) that the set of classes $c\in H^1(X,\ung)$ which have
the same image as $c'$ under the natural map
$$H^1(X,\ung)\to H^1(X,\unh)$$
is equal to the image of the map $i_1$. In particular we deduce
immediately the following claim.
\begin{claim}\label{Claim:001}
If, for any $\ung$-torsor $E'$,
$$H^1(X,\unf(E'))=0$$
then the natural map $H^1(X,\ung)\to H^1(X,\unh)$ is injective.
\end{claim}
We will need the following proposition.

\begin{proposition}\label{P:2.2}
For any $1\leq i<j$ and for any $\ck_0/\ck_j$-torsor $E'$ ,one has
$$H^1(X,(\ck_i/\ck_j)(E'))=0.$$
\end{proposition}
{\bf Proof.} The proof is by induction in $j-i$. Assume first that
$j-i=1$. $(\ck_i/\ck_{i+1})(E')$ is a sheaf of $\co_X$-modules; this
follows easily from (\ref{IsomorSh-1}) and the fact that
$\cw_i/\cw_{i+1}$ is a sheaf of $\co_X$-modules. Hence
$(\ck_i/\ck_{i+1})(E')$ is acyclic.

Assume now that $j-i>1$. Then we have a short exact sequence of
sheaves
$$1\to (\ck_j/\ck_{j-1})(E')\to (\ck_i/\ck_j)(E')\to
(\ck_i/\ck_{j-1})(E')\to 1.$$ Hence we have an exact sequence of
pointed sets (see \cite{groth}, Section 5.3),
$$H^1(X,(\ck_j/\ck_{j-1})(E'))\to H^1(X,(\ck_i/\ck_j)(E'))\to
H^1(X,(\ck_i/\ck_{j-1})(E')).$$ The first and the third terms of the
last sequence vanish by the induction assumption. Hence the middle
term vanishes too. Proposition is proved. \qed

We easily deduce a corollary.
\begin{corollary}\label{C:inject}
The natural map
$$H^1(X,\glv)\to H^1(X,\glo)$$
is injective.
\end{corollary}
{\bf Proof.} By Proposition \ref{P:2.2} $H^1(X,\ck_1(E'))=0$ for any
$\ck_0$-torsor $E'$. Hence, by Claim \ref{Claim:001}, the map
$H^1(X,\ck_0)\to H^1(X,\ck_0/\ck_1)$ is injective. \qed

\hfill

We will need a few more results from \cite{groth}. Assume $X$ is a
paracompact topological space (remind that all our manifolds are
always assumed to be paracompact). Let $\ung$ be a sheaf of groups
on $X$ as before. Let $\unf \triangleleft \ung$ be a subsheaf of
{\itshape normal abelian} subgroups. Let $\unh:=\ung/\unf$ be the
quotient sheaf as before. The action of $\ung$ on $\unf$ by
conjugation induces in this case an action of $\unh$ on $\unf$. For
any $\unh$-torsor $E''$ one has the sheaf $\unf(E'')$ defined
similarly as before. This is a sheaf of abelian groups since $\unf$
is, and $\unh$ acts on $\unf$ by automorphisms. Grothendieck
(\cite{groth}, Section 5.7) has constructed an element
$$\delta E''\in H^2(X,\unf(E''))$$
with the following property: $\delta E''$ vanishes if and only if
the class $[E'']\in H^1(X,\unh)$ lies in the image of the canonical
map $H^1(X,\ung)\to H^1(X,\unh)$.

In order to apply this result in our situation we will need two
lemmas.
\begin{lemma}\label{L:vanish2}
For any $i>0$ and any $\ck_0/\ck_i$-torsor $E''$,
$$H^2(X,(\ck_i/\ck_{i+1})(E''))=0.$$
\end{lemma}
{\bf Proof.} It is easy to see that $(\ck_i/\ck_{i+1})(E'')$ is a
sheaf of $\co_X$-modules. Hence it is acyclic. \qed

\begin{lemma}\label{L:star}
For any $i\geq 0$ the natural map
$$H^1(X,\ck_0/\ck_{i+1})\to H^1(X,\ck_0/\ck_i)$$
is onto.
\end{lemma}
{\bf Proof.} Let $c''\in H^1(X,\ck_0/\ck_i)$ be an arbitrary
element. Let $E''$ be a $\ck_0/\ck_i$-torsor representing $c''$.
Consider the element $\delta E''\in H^2(X,(\ck_i/\ck_{i+1})(E''))$.
Since the last group vanishes by Lemma \ref{L:vanish2}, by the above
mentioned result of Grothendieck, $c''$ lies in the image of
$H^1(X,\ck_0/\ck_{i+1})$. Lemma is proved. \qed

\begin{corollary}\label{C:surject}
The natural map
$$H^1(X,\glv)\to H^1(X,\glo)$$
is onto.
\end{corollary}
{\bf Proof.} The map in the statement factorizes into the sequence
of maps
\begin{eqnarray*}
H^1(X,\glv)=H^1(X,\ck_0)\to H^1(X,\ck_0/\ck_n)\to
H^1(X,\ck_0/\ck_{n-1})\to \dots \\\dots \to
H^1(X,\ck_0/\ck_1)=H^1(X,\glo)
\end{eqnarray*}
where all the maps are surjective by Lemma \ref{L:star}. Hence their
composition is onto too. \qed

Now Theorem \ref{T:biject} follows immediately from Corollaries
\ref{C:inject} and \ref{C:surject}.

\hfill

\begin{remark}
Theorem \ref{T:biject} has the following immediate consequence. For
a compact manifold $X$ we can construct a $K$-ring generated by
finitely generated projective $V^\infty(X)$-modules in the standard
way. Namely as a group it is equal to the quotient of the free
abelian group generated by isomorphism classes of such modules by
the relations
$$[M\oplus N]=[M]+[N].$$
The product is induced by the tensor product of such
$V^\infty(X)$-modules. Then Theorem \ref{T:biject} implies that
there is a canonical isomorphism of this $K$-ring with the classical
topological $K^0$-ring (see \cite{atiyah-lectures}) constructed from
vector bundles.
\end{remark}

\begin{remark}
The main results of this paper are of general nature. It would be
interesting to have concrete geometric examples of
$\cv^\infty_X$-modules; in the classical case of $\co_X$-modules we
have the tangent bundle and its tensor powers.

As a first small step in this direction let us mention the following
construction. Let $\cl$ be a flat vector bundle over a manifold $X$.
By an abuse of notation, we will also denote by $\cl$ the sheaf of
its locally constant sections. Let $\underline{\CC}$ be the constant
sheaf of $\CC$-vector spaces. Consider the $\cv^\infty_X$-module
defined by
$$\tilde\cl:=\cl\otimes_{\underline{\CC}}\cv^\infty_X$$
where we consider $\cv^\infty_X$ as $\underline{\CC}$-module via the
imbedding $\underline{\CC}\inj \cv^\infty_X$ where 1 goes to the
Euler characteristic. It is easy to see that $\tilde\cl$ is a
locally free $\cv^\infty_X$-module.
\end{remark}

\end{document}